\documentclass[12pt]{amsart}
\usepackage{amssymb,latexsym}

\numberwithin{equation}{section}

\newtheorem{theorem}{Theorem}[section]
\newtheorem{proposition}[theorem]{Proposition}

\newtheorem{corollary}[theorem]{Corollary}

\theoremstyle{definition}

\newtheorem{definition}[theorem]{Definition}

\renewcommand{\eqref}[1]{{\rm (\ref{#1})}}

\def\ZZ{\mathbb{Z}}
\def\QQ{\mathbb{Q}}
\def\RR{\mathbb{R}}
\def\B{\mathcal{B}}
\def\N{\mathcal{N}}

\def\endproof{\hfill$\square$\medskip}

\begin{document}

\title[Nested complexes]
{Nested complexes and their polyhedral realizations}

\author{Andrei Zelevinsky}
\address{\noindent Department of Mathematics, Northeastern University,
 Boston, MA 02115}
\email{andrei@neu.edu}

\date{July 18, 2005}

\dedicatory{To Bob MacPherson on the occasion of his 60th birthday}

 \thanks{Research 
supported 
by NSF grant DMS-0500534.}

\subjclass[2000]{Primary 52B}

\maketitle

\section{Introduction}

This note which can be viewed as a complement to \cite{post},
presents a self-contained overview of basic properties of nested complexes
and their two dual polyhedral realizations:
as complete simplicial fans, and as simple polytopes.
Most of the results are not new; our aim is
to bring into focus a striking similarity between nested
complexes and associated fans and polytopes on one side,
and cluster complexes and generalized associahedra
introduced and studied in \cite{yga,cfz} on the other side.

First a very brief history (more details and references can be found
in \cite{tol,post}).
Nested complexes appeared in the work by
De Concini - Procesi \cite{DP} in the context of subspace
arrangements.
More general complexes associated with arbitrary finite graphs
were introduced and studied in \cite{CD} and more recently in
\cite{tol}; these papers also present a construction of
associated simple convex polytopes (dubbed graph associahedra in
\cite{CD}, and De Concini - Procesi associahedra in \cite{tol}).
An even more general setup developed by Feichtner - Yuzvinsky
in \cite{FY} associates
a nested complex and a simplicial fan to an arbitrary finite atomic lattice.
In the present note we follow \cite{FS} and \cite{post}
in adopting an intermediate level of generality, and
studying the nested complexes associated to \emph{building
sets} (see Definition~\ref{def:building} below).

Feichtner and Sturmfels in \cite{FS} show that the simplicial fan
associated to a building set is the normal fan of a simple convex polytope,
while Postnikov in \cite{post} gives
an elegant construction of this polytope as a
Minkowski sum of simplices.
In this note we develop an alternative approach, constructing the
polytope directly from the fan; this construction is completely
analogous to that in \cite{cfz}.
The main role in this approach is played by the link decomposition
property (Proposition~\ref{pr:link-decomposition} below) asserting
that the link of any element in a nested complex can be naturally
identified with a nested complex of smaller rank.
This result provides a unified method for establishing
properties of nested complexes and their polyhedral realizations by induction
on the rank of a building set.

The note is structured as follows.
Section~\ref{sec:nested} introduces (basically following
\cite{post}) building sets and associated nested complexes.
Section~\ref{sec:link-decomposition} presents the link
decomposition property (Proposition~\ref{pr:link-decomposition}).
In this generality it seems to be new; for the graph associahedra
it appeared in \cite{tol}.
Section~\ref{sec:applications} contains some structural properties
of nested complexes, all of them derived in a unified way from
Proposition~\ref{pr:link-decomposition}.

In Section~\ref{sec:fan} we show that any nested complex can be realized
as a smooth complete simplicial fan (Theorem~\ref{th:fan} and
Corollary~\ref{cor:fan}); we refer to this fan as the \emph{nested
fan}.
This construction has appeared in \cite{FY,FK,FS}.
However, it is not so easy to extract from these papers
a proof that the nested fan is well-defined; this is why
we prefer to give a simple self-contained proof.

In Section~\ref{sec:normal-fans}, we prove (Theorem~\ref{th:normal-fan})
that the nested fan is the normal fan of a simple polytope~$\Pi$,
which we call the \emph{nested polytope}.
The argument is completely parallel to that in \cite{cfz}; in
particular, it leads to a specific realization of~$\Pi$ which can
be shown to be the same as the realization given in \cite{post}.

The concluding Section~\ref{sec:graphical} discusses the special class
of \emph{graphical} building sets, those giving rise to graph
associahedra \cite{CD,tol}.
We give a simple characterization of graphical building sets
(Proposition~\ref{pr:graphical-characterization}).
For them an analogy with cluster complexes becomes
sharper: in particular, the corresponding nested complex is
a \emph{clique} complex (Corollary~\ref{cor:clique}).

\section{Nested complexes}
\label{sec:nested}

We start by recalling basic definitions from
\cite[Section~7]{post}.

\begin{definition}[\cite{post}, Definition~7.1]
\label{def:building}
Let~$S$ be a nonempty finite set.
A \emph{building set\/} (or simply a \emph{building\/}) on $S$
is a collection $\B$ of nonempty subsets of $S$
satisfying the conditions:
\begin{itemize}
\item[(B1)]
If $I,J\in \B$ and $I\cap J\neq \emptyset$,
then $I\cup J\in \B$.
\item[(B2)]
$\B$ contains all singletons $\{i\}$, for $i\in S$.
\end{itemize}
\end{definition}

The nested complex associated to a building~$\B$ is defined
as follows.
Let $\B_{\max} \subset \B$ be the set of maximal (by inclusion) elements of $\B$.
We refer to the elements of $\B_{\max}$ as the \emph{$\B$-components}.
In view of (B1) and (B2), the $\B$-components are pairwise disjoint,
and their union is~$S$.
We define the \emph{rank} of~$\B$ by setting
${\rm rk}(\B) = |S| -|\B_{\max}|$.

\begin{definition}[cf.~\cite{post}, Definition~7.3]
\label{def:nested_family}
A subset $N \subset \B - \B_{\max}$ is called \emph{nested}
if it satisfies the following condition:
for any $k \geq 2$ and any $I_1, \dots, I_k \in N$ such that none of the
$I_i$ contains another one, the union $I_1\cup \cdots \cup I_k$ is not in $\B$.
The simplicial complex on $\B - \B_{\max}$ formed by the nested sets
is called the \emph{nested complex\/} and denoted $\N(\B)$.
\end{definition}

In particular, the $1$-simplices of $\N(\B)$ are pairs
$\{I, J\} \subset \B$ such that:
\begin{equation}
\label{eq:B-compatibility}
\text{either $I \subset J$, or
$J \subset I$, or
$I \cup J \notin \B$}
\end{equation}
(in view of (B2), in the last case $I$ and $J$ are disjoint).
Note that our definition of nested sets is \emph{different}
from the one in \cite[Definition~7.3]{post}: Postnikov's nested
sets are obtained from ours by adjoining the set~$\B_{\max}$.
Of course this does not affect the poset structure of the nested complex; an
advantage of the present version is that our nested sets form a simplicial
complex.

\section{Link decomposition}
\label{sec:link-decomposition}

Let~$\B$ be a building on~$S$, and $\N(\B)$ the associated nested complex.
For every $C \in \B - \B_{\max}$, let
\begin{equation}
\label{eq:link}
\N(\B)_C = \{N' \subset \B - \B_{\max} - \{C\} :
N' \cup \{C\} \in \N(\B)\}
\end{equation}
denote the \emph{link} of~$C$ in $\N(\B)$; this is a simplicial
complex on the set $\{I \in \B - \B_{\max} - \{C\} : \{C, I\}
\in \N(\B)\}$.
Our main tool in studying the structure of $\N(\B)$ is the recursive
description of $\N(\B)_C$ given in Proposition~\ref{pr:link-decomposition} below.
To state it, we need some preparation.

\begin{definition}
\label{def:restriction-contraction}
Let~$C$ be a nonempty subset of~$S$.
\begin{itemize}
\item
The \emph{restriction} of~$\B$ to~$C$ is the building on~$C$
defined by
\begin{equation}
\label{eq:restriction}
\B|_C = \{I \subset C : I \in \B\}.
\end{equation}
\item
The \emph{contraction} of~$C$ from~$\B$ is the building on~$S - C$
defined by
\begin{equation}
\label{eq:contraction}
C \backslash \B = \{I \subset S - C : I \in \B \,\, \text{or} \,\,
C \cup I \in \B\}.
\end{equation}
\end{itemize}
\end{definition}

The fact that both $\B|_C$ and $C \backslash \B$ are indeed
buildings, is immediate from Definition~\ref{def:building}.

Let $S_1, \dots, S_k$ be disjoint finite sets, and let $\B_i$ be a
building on~$S_i$ for $i = 1, \dots, k$.
The \emph{product} $\B_1 \times \cdots \times \B_k$ is defined as the
building on $S_1 \sqcup \cdots \sqcup S_k$ formed by the elements
of all the~$\B_i$.
Clearly, the nested complex $\N(\B_1 \times \cdots \times \B_k)$ is
naturally isomorphic to the product
$\N(\B_1) \times \cdots \times \N(\B_k)$.
In particular, every building~$\B$ on~$S$ decomposes as a product
$\B|_{S_1} \times \cdots \times \B|_{S_k}$, where $S_1, \dots, S_k$ are the
$\B$-components; therefore,
$\N(\B)$ can be identified with $\N(\B|_{S_1}) \times \cdots \times \N(\B|_{S_k})$.

\begin{proposition}[Link Decomposition]
\label{pr:link-decomposition}
For every $C \in \B - \B_{\max}$, the link $\N(\B)_C$ is
isomorphic to $\N(\B|_C \times C \backslash \B)$.
\end{proposition}

\proof
Fix $C \in \B - \B_{\max}$ and abbreviate
$\B' = \B|_C \times C \backslash \B$.
Let
$$\B_0 = \{I \in \B - \B_{\max} - \{C\} : \{C, I\} \in \N(\B)\}$$
be the ground set of $\N(\B)_C$.
In view of \eqref{eq:B-compatibility}, $\B_0$
is the disjoint union of the
following three sets:
\begin{align}
\label{eq:B1}
& \B_1 = \{I \in \B - \{C\} : I \subset C\};\\
\label{eq:B2}
& \B_2 = \{I \in \B - \B_{\max} - \{C\} : C \subset I\};\\
\label{eq:B3}
& \B_3 = \{I \in \B - \B_{\max} : C \cap I = \emptyset, \,\, C \cup I \notin \B\}.
\end{align}
Now to every $I \in \B_0 = \B_1 \cup \B_2 \cup \B_3$ we associate
a subset $I' \subset S$ by setting
\begin{equation}
\label{eq:ItoI'}
I' =
\begin{cases}
I & \text{if $I \in \B_1 \cup \B_3$;} \\[.05in]
I - C & \text{if $I \in \B_2$.}
\end{cases}
\end{equation}
Remembering the definitions, we see that $\B_1$ is the ground set
of the nested complex $\N(\B|_C)$, and $\B_3$ is a subset of
the ground set of $\N(C \backslash \B)$ formed by those
$I' \subset S - C$ for which $C \cup I' \notin \B$.
Furthermore, the correspondence $I \mapsto I' = I - C$ identifies
$\B_2$ with the remaining part of the ground set of $\N(C \backslash \B)$,
formed by those $I' \subset S - C$ for which $C \cup I' \in \B$.
Thus, the correspondence $I \mapsto I'$ is a bijection between the
ground sets of~$\N(\B)_C$ and~$\N(\B')$.
To prove Proposition~\ref{pr:link-decomposition}, it remains to show the following:
\begin{eqnarray}
\label{eq:ItoI'-isomorphism}
& \text{The correspondence $I \mapsto I'$ induces an}\\
\nonumber
& \text{isomorphism of complexes $\N(\B)_C$ and $\N(\B')$.}
\end{eqnarray}

Thus, we need to show the following: for any distinct
$I_1, \dots, I_k \in \B_0$, the set $\{C, I_1, \dots, I_k\}$
is nested for~$\B$ if and only if $\{I'_1, \dots, I'_k\}$
is nested for~$\B'$.
First suppose $\{I'_1, \dots, I'_k\}$ is \emph{not}
nested for~$\B'$.
Without loss of generality, we can assume that $k \geq 2$,
none of the $I'_i$ contains another one, and
$I'_1 \cup \cdots \cup I'_k \in \B'$.
Choosing the smallest possible~$k$ with this property and
using property (B1) in Definition~\ref{def:building}, we can
further assume that the $I'_i$ are pairwise disjoint.
Remembering the definition of~$\B'$, we see that either all
the $I'_i$ are contained in~$C$, or all
the $I'_i$ are contained in~$S - C$.
In the former case, $I_i = I'_i$ for all~$i$, so
we have $I_1 \cup \cdots \cup I_k \in \B$,
implying that $\{C, I_1, \dots, I_k\}$
is not nested for~$\B$.
Now suppose that all the $I'_i$ are contained in~$S - C$.
Remembering \eqref{eq:contraction}, we see that the condition
$I'_1 \cup \cdots \cup I'_k \in \B'$ means that at least one of
the subsets $I'_1 \cup \cdots \cup I'_k$ and
$C \cup I'_1 \cup \cdots \cup I'_k$ belongs to~$\B$.
If~$I'_i = I_i$ for all~$i = 1, \dots, k$, we again conclude that
$\{C, I_1, \dots, I_k\}$ is not nested for~$\B$.
Therefore suppose that $I_i = C \cup I'_i$ for some~$i$.
If $I_j = C \cup I'_j$ for some $i \neq j$,
then $I_i \cup I_j \in \B$ by
property (B1) in Definition~\ref{def:building}, and so
$\{C, I_1, \dots, I_k\}$ is not nested for~$\B$ in this case as well.
Finally, if $I_j = I'_j$
for $j \neq i$, then $I_1, \dots, I_k$ are pairwise disjoint, and
their union belongs to~$\B$, again implying that
$\{C, I_1, \dots, I_k\}$ is not nested for~$\B$.

The converse implication (if $\{C, I_1, \dots, I_k\}$
is not nested for~$\B$ then $\{I'_1, \dots, I'_k\}$
is not nested for~$\B'$) follows by reversing the above arguments.
\endproof

\section{Some applications of the link decomposition}
\label{sec:applications}

\begin{proposition}
\label{pr:purity}
The nested complex $\N(\B)$ is pure of dimension
${\rm rk}(\B) - 1$, that is, all maximal nested sets are of
the same cardinality ${\rm rk}(\B)$.
\end{proposition}

\proof
Proceed by induction on ${\rm rk}(\B)$.
The statement is trivial if ${\rm rk}(\B) = 0$:
in this case, $\B$ consists of the singletons, and $\N(\B)$ is just empty.
Now assume that ${\rm rk}(\B) \geq 1$, and let~$N \subset \B - \B_{\max}$
be a maximal nested set.
Let $C \in N$.
Inspecting Definition~\ref{def:restriction-contraction}, we see
that $|(\B|_C \times C \backslash \B)_{\max}| = |\B_{\max}| + 1$
(the $\B$-components that are disjoint from~$C$ remain
intact, while the component~$S'$ containing~$C$ splits into the
two $(\B|_C \times C \backslash \B)$-components~$C$ and $S' - C$).
Thus, ${\rm rk}(\B|_C \times C \backslash \B) = {\rm rk}(\B) - 1$.
Using Proposition~\ref{pr:link-decomposition} and the induction assumption,
we conclude that
$$|N| = 1 + {\rm rk}(\B|_C \times C \backslash \B) = {\rm rk}(\B),$$
 as claimed.
\endproof

\begin{proposition}
\label{pr:exchanges}
For every maximal nested set~$N$ and every $I \in N$,
there is precisely one maximal nested set~$N'$ such that
$N \cap N' = N - \{I\}$.
\end{proposition}

\proof
First assume that ${\rm rk}(\B) = 1$.
Then $\B$ consists of all
singletons in~$S$ and just one two-element set $\{i,j\}$.
Thus, the only two maximal nested sets are $N = \{I\}$ and
$N' = \{J\}$, where $I$ (resp.~$J$) is the singleton $\{i\}$
(resp.$\{j\}$).
So our assertion holds.

Now let ${\rm rk}(\B) > 1$.
Choose any $C \in N - \{I\}$.
To prove Proposition~\ref{pr:exchanges}, it suffices to show that
the assertion holds for the link $\N(\B)_C$.
It remains to apply Proposition~\ref{pr:link-decomposition} and
induction on ${\rm rk}(\B)$.
\endproof

\begin{definition}
\label{def:dual-graph}
The \emph{dual graph} of the nested complex $\N(\B)$ has maximal
nested sets as vertices, with~$N$ and~$N'$ joined by an edge
whenever $|N \cap N'| = |N| - 1$.
\end{definition}

The following property is immediate from
Proposition~\ref{pr:exchanges}.

\begin{proposition}
\label{pr:dual-graph-regular}
The dual graph of $\N(\B)$ is regular of degree ${\rm rk}(\B)$.
\end{proposition}

We now give a more detailed description of the edges of the dual
graph.

\begin{proposition}
\label{pr:dual-graph-exchanges}
Let $N_1$ and $N_2$ be two
maximal nested sets in $\N(\B)$ joined by an edge in the dual
graph, and let $N_1 \cap N_2 = N_1 - \{I_1\} = N_2 - \{I_2\}$.
Then the following properties hold:
\begin{enumerate}
\item
Neither of $I_1$ and $I_2$ contains another one.
\item
If $I_1 \cap I_2 \neq \emptyset$ then all the
$\B|_{I_1 \cap I_2}$-components belong to $N_1 \cap N_2$.
\item
There exist pairwise disjoint $I_3, \dots, I_k \in N_1 \cap N_2$
such that $(I_1 \cup I_2) \cap I_i = \emptyset$ for $3 \leq i \leq k$,
and $I_1 \cup \cdots \cup I_k \in (N_1 \cap N_2) \cup \B_{\max}$
(note that
the family $\{I_3, \dots, I_k\}$ can be empty).
\end{enumerate}
\end{proposition}

\proof
As before, we proceed by induction on ${\rm rk}(\B)$.
If ${\rm rk}(\B) = 1$ then $I_1 = \{i\}$ and $I_2 = \{j\}$, with
$\{i, j\}$ the only non-singleton member of~$\B$.
This makes the properties $(1) - (3)$ obvious.

So we assume that ${\rm rk}(\B) > 1$, choose $C \in N_1 \cap N_2$,
and assume that $(1) - (3)$ hold for the building
$\B' = \B|_C \times C \backslash \B$ and its maximal nested sets
$N'_1$ and $N'_2$ corresponding to $N_1$ and $N_2$ as
in Proposition~\ref{pr:link-decomposition}.
Thus, for $i \in \{1, 2\}$, we have
$N'_i = \{I' : I \in N_i - \{C\}\}$, where the correspondence
$I \mapsto I'$ is given by \eqref{eq:ItoI'}.
For the sake of convenience, we refer to the properties
$(1) - (3)$ for the sets $N'_1$ and $N'_2$ as $(1') - (3')$.
So we need to show that $(1') - (3')$ imply $(1) - (3)$.

First of all, by the definition of $\B'$, every $J \in \B'$ is
contained in either~$C$ or $S - C$.
Hence $(3')$ implies that either both~$I'_1$ and~$I'_2$
are contained in~$C$ or both are contained in~$S - C$.
Remembering \eqref{eq:ItoI'}, we see that $(1')$ implies (1).

To prove (2) suppose that $J = I_1 \cap I_2 \neq \emptyset$.
There are three possibilities to consider:
\begin{itemize}
\item
$J$ is strictly contained in~$C$.
\item
$J$ contains~$C$.
\item
$J$ is disjoint from~$C$.
\end{itemize}

If~$J$ is strictly contained in~$C$ then both $I_1$ and $I_2$
belong to the set $\B_1$ in \eqref{eq:B1}.
By \eqref{eq:ItoI'}, we have $I'_i = I_i$ for $i \in \{1, 2\}$,
hence $I'_1 \cap I'_2 = J$.
Again by \eqref{eq:ItoI'}, $\B|_J = \B'|_j$, and so
$(2')$ implies (2).

If~$J$ contains~$C$ then both $I_1$ and $I_2$
belong to the set $\B_2$ in \eqref{eq:B2}.
By \eqref{eq:ItoI'}, we have $I'_i = I_i - C$ for $i \in \{1, 2\}$,
hence $J' = I'_1 \cap I'_2 = J - C$.
If $J' = \emptyset$ then $J = C$ is the only $\B|_J$-component,
making (2) trivial.
So suppose $J'$ is nonempty, and let $J_1, \dots, J_p$ be the
$\B'|_{J'}$-components.
By property (B1) in Definition~\ref{def:building},
there can be at most one of
the $J_i$ such that $C \cup J_i \in \B$.
If $C \cup J_i \notin \B$ for all~$i$, then the
$\B|_J$-components are $C, J_1, \dots, J_p$;
if say $C \cup J_1 \in \B$, then the
$\B|_J$-components are $C \cup J_1, J_2, \dots, J_p$.
In both cases, (2) follows from $(2')$, as desired.

It remains to prove (3).
Since by Proposition~\ref{pr:link-decomposition}, every member of
$N'_1 \cap N'_2$ is of the form~$I'$ for
$I \in (N_1 \cap N_2) - \{C\}$, we assume that $(3')$ holds with the
subsets $I'_3, \dots, I'_k \in N'_1 \cap N'_2$ in place of
$I_3, \dots, I_k \in (N_1 \cap N_2) - \{C\}$.

By the definition of~$\B'$, $(3')$ implies that $I'_1, \dots, I'_k$
either are all contained in~$C$, or all contained in~$S - C$.
In the former case, $I'_i = I_i$ for $i = 1, \dots, k$, hence
(3) follows at once from $(3')$.
So suppose that $I'_i \subset S - C$ for $i = 1, \dots, k$.

Since by $(3')$ the sets $I'_1, I'_3, \dots, I'_k$ belong to~$N'_1$,
they form a nested set, hence $C \cup I'_i \in \B$ for at most one
index $i \in \{1, 3, \dots, k\}$.
By the same token, $C \cup I'_i \in \B$ for at most one
index $i \in \{2, 3, \dots, k\}$.
In other words, $I_i = C \cup I'_i$ for at most one index~$i$
among $\{1, 3, \dots, k\}$, and at most one among $\{2, 3, \dots, k\}$,
while $I_i = I'_i$ for the rest of the indices.
Therefore the conditions that $I_3, \dots, I_k$ are
pairwise disjoint, and that $(I_1 \cup I_2) \cap I_i = \emptyset$
for $3 \leq i \leq k$, follow from the corresponding conditions
for the~$I'_i$.

To finish the proof, it remains to consider the following
two cases:
\begin{itemize}
\item
$I_i = I'_i$ for $i = 1, \dots, k$.
\item
$I_i = C \cup I'_i$ for some $i \in \{1, \dots, k\}$.
\end{itemize}

First suppose $I_i = I'_i$ for $i = 1, \dots, k$.
Then the  condition that
$I'_1 \cup \cdots \cup I'_k \in (N'_1 \cap N'_2) \cup \B'_{\max}$
means that either $I_1 \cup \cdots \cup I_k \in (N_1 \cap N_2) \cup \B_{\max}$,
or $C \cup I_1 \cup \cdots \cup I_k \in (N_1 \cap N_2) \cup \B_{\max}$.
In the former (resp. latter) case, (3) holds for $I_3, \dots, I_k$
(resp. for $I_3, \dots, I_k, I_{k+1} = C$).

Finally, if $I_i = C \cup I'_i$ for some $i \in \{1, \dots, k\}$ then
$I_1 \cup \cdots \cup I_k = C \cup I'_1 \cup \cdots \cup I'_k \in \B$
by property (B1) in Definition~\ref{def:building}.
Therefore, $I'_1 \cup \cdots \cup I'_k =
(I_1 \cup \cdots \cup I_k)'$, hence
$I_1 \cup \cdots \cup I_k \in (N_1 \cap N_2) \cup \B_{\max}$,
finishing the proof.
\endproof

Let~$D$ be a nested set of cardinality ${\rm rk}(\B) - 2$.
As in~\cite[Section~2.1]{ca2}, we call the induced subgraph
of the dual graph whose vertices are the maximal nested
sets containing~$D$, the \emph{geodesic loop} associated to~$D$.
The geodesic loop can be identified with the dual graph of the
link of~$D$ in $\N(\B)$.
By a repeated application of Proposition~\ref{pr:link-decomposition},
a geodesic loop is the dual graph of some complex $\N(\B')$
with ${\rm rk}(\B') = 2$.
Thus, by Proposition~\ref{pr:dual-graph-regular}, any geodesic
loop is a $2$-regular graph, hence a disjoint union of cycles.
More precisely, we have the following opportunities.

\begin{proposition}
\label{pr:geodesic-loops}
Any geodesic loop is a $d$-cycle for some $d \in \{3, 4, 5, 6\}$.
\end{proposition}

\proof
We can assume that ${\rm rk}(\B) = 2$, hence
a geodesic loop in question is the dual graph of $\N(\B)$.
Without affecting this graph, we can assume that there are no
singleton $\B$-components of~$S$.
An easy inspection then leaves us with the following options
for~$\B$.

\begin{enumerate}
\item[(D1)]
$S = \{1, 2, 3\}$; $\B$ consists of the singletons and~$S$.
The dual graph is a $3$-cycle with vertices (= maximal nested sets)
$$\{\{1\},\, \{2\}\},\, \{\{2\}, \{3\}\},\, \{\{1\}, \{3\}\}.$$
\item[(D2)]
$S = \{1, 2, 3, 4\}$; $\B$ consists of the singletons,
$\{1,3\}$ and $\{2,4\}$.
The dual graph is a $4$-cycle with vertices
$$\{\{1\}, \{2\}\},\, \{\{2\}, \{3\}\},\,
\{\{3\}, \{4\}\},\, \{\{1\}, \{4\}\}.$$
\item[(D3)]
$S = \{1, 2, 3\}$; $\B$ consists of the singletons, $\{1,2\}$, and~$S$.
The dual graph is a $4$-cycle with vertices
$$\{\{1\}, \{1,2\}\},\, \{\{2\}, \{1,2\}\},\, \{\{2\}, \{3\}\},\,
\{\{1\}, \{3\}\}.$$
\item[(D4)]
$S = \{1, 2, 3\}$; $\B$ consists of the singletons, $\{1,2\}$, $\{2,3\}$, and~$S$.
The dual graph is a $5$-cycle with vertices
$$\{\{1\}, \{1,2\}\},\, \{\{2\}, \{1,2\}\},\,
\{\{2\}, \{2, 3\}\},\, \{\{3\}, \{2,3\}\}, \, \{\{1\}, \{3\}\}.$$
\item[(D5)]
$S = \{1, 2, 3\}$; $\B$ consists of all nonempty subsets of~$S$.
The dual graph is a $6$-cycle with vertices
$$\{\{1\}, \{1,2\}\},\, \{\{2\}, \{1,2\}\},\,
\{\{2\}, \{2,3\}\},\, \{\{3\}, \{2,3\}\},\, \{\{3\},
\{1,3\}\},\, \{\{1\}, \{1,3\}\}.$$
\end{enumerate}
\endproof

We conclude this section by using Proposition~\ref{pr:link-decomposition}
to obtain a recursion for the $f$-vector of a nested complex.
For $k \geq 0$, let $f_k (\N(\B))$ denote the number of
nested $k$-subsets of $\N(\B)$ (so that, in particular,
$f_0 (\N(\B)) = 1$, and $f_1 (\N(\B)) = |\B| - |\B_{\max}|$).
Let
$$f(\N(\B)) = \sum_{k \geq 0} f_k (\N(\B)) x^k$$
be the corresponding generating function.
The following result is an analogue of
\cite[Proposition~3.7]{yga}.

\begin{proposition}
\label{pr:f-vector}
\begin{enumerate}
\item
We have
$f (\N(\B_1 \times \B_2)) = f(\N(\B_1)) f(\N(\B_2))$.

\item
For every $k \geq 1$, we have
\begin{equation}
\label{eq:f-number}
k f_k (\N(\B)) =
\sum_{C \in \B - \B_{\max}}
f_{k-1}(\N(\B|_C \times C \backslash \B)) \,.
\end{equation}
Equivalently,
$$\frac{d f(\N(\B))}{d x} =
\sum_{C \in \B - \B_{\max}} f(\N(\B|_C))f(\N(C \backslash \B)) \ .$$
\end{enumerate}
\end{proposition}

\proof
Part (1) is immediate from the definitions.
To prove (2) count in two different ways the number of pairs
$(C, N)$, where $C \in \B - \B_{\max}$, and~$N$ is a nested
$k$-set containing~$C$.
On the one hand, each~$N$ contributes~$k$ to this count, so the number in question
is equal to $k f_k (\N(\B))$; on the
other hand, by Proposition~\ref{pr:link-decomposition}, each~$C$ contributes
$f_{k-1}(\N(\B|_C \times C \backslash \B))$, adding up to the right hand
side of \eqref{eq:f-number}.
\endproof

It would be interesting to compare the recursion for the
$f$-vector given by \eqref{eq:f-number} with the one
in \cite[Theorem~7.11 (3)]{post}.

\section{Nested fans}
\label{sec:fan}

We now realize the nested complex $\N(\B)$ associated to a
building~$\B$ on a finite set~$S$ as a smooth complete simplicial
fan in a real vector space of dimension ${\rm rk}(\B)$.
Let $\RR^S$ (resp.~$\ZZ^S$) be a real vector space (resp. a
lattice) with the distinguished basis $\{e_i : i \in S\}$.
For a subset $I \subset S$, we denote
$$e_I = \sum_{i \in I} e_i \in \ZZ^S \subset \RR^S.$$
Let~$V$ denote the quotient space of~$\RR^S$ modulo the linear span
of the vectors $e_C$ for $C \in \B_{\max}$ (recall that $\B_{\max}$
stands for the set of all $\B$-components of~$S$, that is, the
maximal elements of~$\B$).
Thus, $\dim V = {\rm rk}(\B)$.
Let $\pi: \RR^S \to V$ denote the projection, and let
$L = \pi(\ZZ^S) \subset V$ be the standard lattice in~$V$.
We abbreviate $\overline e_I = \pi(e_I) \in L$.

\begin{theorem}
\label{th:fan}
\begin{enumerate}
\item
For every maximal nested set~$N$, the vectors
$\overline e_I$ for $I \in N$ form a $\ZZ$-basis in~$L$.
\item
Every vector~$v \in L$ has a unique expansion (referred to
as the \emph{nested expansion})
$v = \sum_{I \in \B - \B_{\max}} c_I \overline e_I$ such that all
coefficients~$c_I$ are nonnegative integers, and the set
$\{I : c_I > 0\}$ is nested.
\end{enumerate}
\end{theorem}

\proof
(1) It suffices to show that any maximal nested set~$N$ satisfies:
\begin{equation}
\label{eq:basis}
\text{The vectors $e_I$ for $I \in N \cup \B_{\max}$
form a $\ZZ$-basis in~$\ZZ^S$.}
\end{equation}
We proceed by induction on $|N| = {\rm rk}(\B)$
(see Proposition~\ref{pr:purity}).
If ${\rm rk}(\B) = 0$ then~$N$ is empty, and
$\B = \B_{\max}$ consists of the singletons, making
\eqref{eq:basis} trivial.
If ${\rm rk}(\B) > 0$, choose any $C \in N$ and apply the
induction assumption to the building $\B' = \B|_C \times C \backslash \B$
on~$S$ and the maximal nested set $N' \in \N(\B|_C \times C \backslash \B)$
corresponding to~$N$ as in Proposition~\ref{pr:link-decomposition}.
Replacing $(\B, N)$ with $(\B', N')$ changes the vectors in \eqref{eq:basis}
as follows: for each $I \in N \cup \B_{\max}$ that strictly
contains~$C$, replace $e_I$ with $e_{I - C} = e_I - e_C$.
Clearly, this transformation and its inverse preserve the set of
$\ZZ$-bases in~$\ZZ^S$, and we are done.

(2) Any~$v \in L$ can be uniquely expressed as
$v = \pi(\sum_{i \in S} c_i e_i)$, where all $c_i$ are nonnegative
integers, and $\min_{i \in C} c_i = 0$ for every
$C \in \B_{\max}$.
Let $S_1 \supset S_2 \supset \cdots$ be the decreasing sequence of
subsets of~$S$ given by $S_j = \{i \in S : c_i \geq j\}$.
Then we have
\begin{equation}
\label{eq:expansion}
v = \sum_{j \geq 1} \sum_{I \in (\B|_{S_j})_{\max}}
\overline e_I.
\end{equation}
Collecting similar terms in \eqref{eq:expansion} yields an
expansion of~$v$ into a nonnegative integer linear combination of
the~$\overline e_I$ for $I \in \B - \B_{\max}$.
It is easy to see that this expansion is nested, i.e.,
the contributing sets~$I$ form a nested set.

Now let $v = \sum_{I \in \B - \B_{\max}} c_I \overline e_I$ be a
nested expansion, and let $N = \{I : c_I > 0\}$ be the
corresponding nested set.
Let $I_1, \dots, I_k$ be the maximal (by inclusion) elements
of~$N$, and let $S_1 = I_1 \cup \cdots \cup I_k$.
As an easy consequence of the definition of a nested set,
the sets $I_1, \dots, I_k$ are the $\B|_{S_1}$-components; in
particular, we see that $S_1$ does not contain any
$\B$-component.
Applying the same construction to the vector
$v' = v - \sum_{j = 1}^k \overline e_{I_j}$, and continuing in the
same manner, we conclude that the nested expansion of~$v$ that we started
with, is given by \eqref{eq:expansion}.
This proves the uniqueness of a nested expansion, and concludes the proof of
Theorem~\ref{th:fan}.
\endproof

For every nested set $N \in \N(\B)$, let $\RR_{\geq 0} N$
denote the cone in~$V$ generated by the vectors
$\overline e_I$ for $I \in N$.
Theorem~\ref{th:fan}
has the following geometric corollary.

\begin{corollary}
\label{cor:fan}
The cones $\RR_{\geq 0} N$ for $N \in \N(\B)$ form a smooth
complete simplicial fan in~$V$.
\end{corollary}

We call the fan in Corollary~\ref{cor:fan} the \emph{nested fan}
associated to~$\B$, and denote it by~$\Delta(\B)$.
Corollary~\ref{cor:fan} follows from Theorem~\ref{th:fan} by
fairly standard arguments (see e.g.,
\cite[Section~3.4, proof of Theorem~1.10]{yga}).
For the reader's convenience, here is a self-contained proof.

\smallskip

\noindent {\it Proof of Corollary~\ref{cor:fan}.}
By Theorem~\ref{th:fan} (1), every cone $\RR_{\geq 0} N$ for $N \in \N(\B)$
is simplicial and is generated by a part of a $\ZZ$-basis in the
lattice~$L$; in other words, if these cones form a fan, then
this fan is \emph{simplicial} and \emph{smooth}.
To check that the nested fan is well defined, we need to show that
\begin{equation}
\label{eq:cone-intersection}
\RR_{\geq 0} N_1 \cap \RR_{\geq 0} N_2 = \RR_{\geq 0} (N_1 \cap N_2)
\quad (N_1, N_2 \in \N(\B)) \ ;
\end{equation}
the completeness property asserts that
\begin{equation}
\label{eq:cone-union}
\bigcup_{N \in \N(\B)} \RR_{\geq 0} N = V \ .
\end{equation}

Let us first show \eqref{eq:cone-union}.
The existence of a nested expansion in Theorem~\ref{th:fan} (2)
means that $\bigcup_{N \in \N(\B)} \ZZ_{\geq 0} N = L$, and so
the union $\bigcup_{N \in \N(\B)} \QQ_{\geq 0} N$ is the
$\QQ$-vector space generated by~$L$.
Therefore, the union $\bigcup_{N \in \N(\B)} \RR_{\geq 0} N$
is dense in~$V$; since it is also closed, \eqref{eq:cone-union}
follows.

Turning to \eqref{eq:cone-intersection}, we need to show that
$$\RR_{\geq 0} N_1 \cap \RR_{\geq 0} N_2 \subset \RR_{\geq 0} (N_1 \cap N_2)$$
(the reverse inclusion is trivial).
Clearly, it suffices to show the following:
\begin{equation}
\label{eq:RR-uniqueness}
\text{If $N_1, N_2 \in \N(\B)$ and
$\RR_{> 0} N_1 \cap \RR_{> 0} N_2 \neq \emptyset$, then $N_1 = N_2$,}
\end{equation}
where $\RR_{> 0} N$ stands for the set of all \emph{positive}
linear combinations of the vectors $\overline e_I$ for $I \in N$.
Let $\RR N$ denote the linear span of the cone $\RR_{\geq 0} N$.
Since each of the cones $\RR_{\geq 0} N_1$ and $\RR_{\geq 0} N_2$
is generated by a part of a $\ZZ$-basis in~$L$, the vector
subspace $\RR N_1 \cap \RR N_2 \subset V$ is defined over~$\ZZ$.
The intersection $\RR_{> 0} N_1 \cap \RR_{> 0} N_2$ is an open polyhedral cone
in $\RR N_1 \cap \RR N_2$.
Therefore, if $\RR_{> 0} N_1 \cap \RR_{> 0} N_2 \neq \emptyset$
then $\ZZ_{> 0} N_1 \cap \ZZ_{> 0} N_2 \neq \emptyset$, and so
$N_1 = N_2$ by the uniqueness of a nested expansion in
Theorem~\ref{th:fan} (2).
\endproof

As a standard consequence of Corollary~\ref{cor:fan}, we get
the following corollary.

\begin{corollary}
\label{cor:dual-graph-connected}
The dual graph (see Definition~\ref{def:dual-graph}) of any nested
complex is connected.
\end{corollary}

\section{Nested polytopes}
\label{sec:normal-fans}

\begin{theorem}
\label{th:normal-fan}
The nested fan $\Delta(\B)$ in Corollary~\ref{cor:fan} is the normal fan
of a simple convex polytope~$\Pi(\B)$,
which can be explicitly realized 
as the set of all tuples $(x_i)_{i \in S}$ of real numbers
satisfying:
\begin{align}
\label{eq:equalities}
\sum_{i \in C} x_i & \, = \, 0 \quad (C \in \B_{\max}) \, ,\\
\label{eq:inequalities}
\sum_{i \in I} x_i &\,  \leq \, |I|\ (2^{|C|-1} - 2^{|I|-1}) \quad
(I \in \B - \B_{\max}, \,\, I \subset C \in \B_{\max}) \, .
\end{align}
\end{theorem}

Before proving Theorem~\ref{th:normal-fan}, we recall relevant facts
on normal fans (see \cite{ziegler} or
\cite[Section~2.1]{cfz}).

\begin{definition}
\label{def:normal-fan}
Let~$\Pi$ be a full-dimensional simple convex polytope in a real vector
space~$V^*$.
The \emph{support function} of $\Pi$ is
a real-valued function $F$ on the dual vector space~$V$ given by
$$F(v) = \max_{\varphi \in \Pi} \langle v, \varphi\rangle \ .$$
The \emph{normal fan} 
of~$\Pi$ is the complete simplicial fan in~$V$
whose 
cones are in the following inclusion-reversing bijection with the faces
of~$\Pi$:
each face $\Sigma \subset \Pi$ gives rise to the cone
\[\mathcal{C}_\Sigma =
\{v \in V : \langle v,\varphi\rangle = F(v)
\,\, \text{for} \,\, \varphi \in \Sigma\} \, ,
\]
whose codimension is equal to~$\dim (\Sigma)$.
\end{definition}

Now let $\Delta$ be a complete simplicial fan in $V$,
and let~$R$ be a set of
representatives of $1$-dimensional cones in~$\Delta$.
The following proposition (cf.~\cite[Lemma~2.1]{cfz})
gives a criterion for~$\Delta$ to be the normal fan of a simple
convex polytope~$\Pi \subset V^*$.

\begin{proposition}
\label{pr:normal-fan-criterion}
The following conditions are equivalent:
\begin{enumerate}
\item
$\Delta$ is the normal fan 
of a full-dimensional simple convex polytope~$\Pi$ in $V^*$.
\item
There exists a real-valued function~$F$ on~$R$
satisfying the following
system of linear inequalities.
For each pair of adjacent maximal cones $\mathcal {C}$ and $\mathcal
{C}'$ in $\Delta$,
let $\{e\} = (R \cap \mathcal {C}) - (R \cap \mathcal {C}')$ and
$\{e'\} = (R \cap \mathcal {C}') - (R \cap \mathcal {C})$, and
write the unique (up to a nonzero real multiple) linear
dependence between the elements of
$R \cap (\mathcal {C} \cup \mathcal {C}')$ in the form
\begin{equation}
\label{eq:dependence}
m_e e + m_{e'} e' -
\sum_{v \in R \cap \mathcal {C} \cap \mathcal {C}'}
m_v v = 0 \, ,
\end{equation}
where $m_e$ and $m_{e'}$ are positive real numbers.
Then
\begin{equation}
\label{eq:dependence-convexity}
m_e F(e) + m_{e'} F(e') -
\sum_{v \in R \cap \mathcal {C} \cap \mathcal {C}'}
m_v F(v) > 0 \, .
\end{equation}
\end{enumerate}
Under these conditions, one can choose the polytope~$\Pi$ to be defined by the
following system of linear inequalities:
\begin{equation}
\label{eq:Pi-general}
\Pi = \{\varphi \in V^* : \langle e, \varphi \rangle \leq
F(e) \,\, \text{for} \,\, e \in R\} \ .
\end{equation}
\end{proposition}

\noindent {\it Proof of Theorem~\ref{th:normal-fan}.}
We apply the criterion in Proposition~\ref{pr:normal-fan-criterion}
to the nested fan $\Delta(\B)$.
We retain the notation in Section~\ref{sec:fan}.
Thus, the ambient space~$V$ of $\Delta(\B)$ is
the quotient space of~$\RR^S$ modulo the linear span
of the vectors $e_C$ for $C \in \B_{\max}$.
Therefore, the dual space~$V^*$ is naturally identified
with the space of real tuples $(x_i)_{i \in S}$
satisfying \eqref{eq:equalities}.
It remains to show that the linear inequalities
\eqref{eq:inequalities} can be identified
with those in \eqref{eq:Pi-general},
for an appropriate choice of the set of representatives~$R$,
and the function~$F$ on~$R$ satisfying
\eqref{eq:dependence-convexity}.

By the definition of the nested fan, the set~$R \subset V$
in Proposition~\ref{pr:normal-fan-criterion}
can be chosen as follows:
$R = \{\overline e_I : I \in \B - \B_{\max}\}$.
Applying Proposition~\ref{pr:dual-graph-exchanges} and using the notation
introduced there, we see that every relation of the
form \eqref{eq:dependence} can be written as follows:
\begin{equation}
\label{eq:dependence-concrete}
\overline e_{I_1} + \overline e_{I_2} -
\sum_{J \in (\B|_{I_1 \cap I_2})_{\max}}
\overline e_{J} + \sum_{i=3}^k \overline e_{I_i}
- \overline e_{I_1 \cup \cdots \cup I_k} = 0 \ .
\end{equation}
The corresponding inequality \eqref{eq:dependence-convexity}
takes the form
\begin{equation}
\label{eq:convexity-concrete}
\sum_{i=1}^k F(\overline e_{I_i}) -
\sum_{J \in (\B|_{I_1 \cap I_2})_{\max}}
F(\overline e_{J})
- F(\overline e_{I_1 \cup \cdots \cup I_k}) > 0 \ .
\end{equation}
Note that all the sets participating in \eqref{eq:convexity-concrete}
are contained in $I_1 \cup \cdots \cup I_k \in \B$, hence they are
contained in the same $\B$-component~$C$.
Thus, to prove Theorem~\ref{th:normal-fan}, it suffices to show
that \eqref{eq:convexity-concrete} holds for
the function~$F$ given by
$$F(\overline e_{I}) =  |I|\ (2^{|C|-1} - 2^{|I|-1})$$
for all participating subsets~$I$.
In other words, using the
abbreviation
$$f_n(m) = m(2^{n-1} - 2^{m-1}) \ ,$$
it is enough to show that, for any positive integer~$n$,
this function  satisfies
\begin{equation}
\label{eq:f-convexity}
\sum_{i=1}^k f_n(p_i) - \sum_{j=1}^\ell f_n(q_j) - f_n(r) > 0
\end{equation}
for all nonnegative integers $p_i, q_j$, and $r$ not exceeding~$n$ and such that
\begin{equation}
\label{eq:f-conditions}
\sum_{i=1}^k p_i - \sum_{j=1}^\ell q_j - r = 0,
\quad \sum_{j=1}^\ell q_j  < \min(p_1, p_2) \ .
\end{equation}

To prove \eqref{eq:f-convexity}, note that, for every positive
integers~$a$ and~$b$, we have
\begin{equation}
\label{eq:f-convexity-2}
f_n(a) + f_n(b) - f_n(a+b) = a \ (2^{a+b-1} -  2^{a-1}) +
b \ (2^{a+b-1} -  2^{b-1}) > 0 \ .
\end{equation}
It follows that, under the assumptions \eqref{eq:f-conditions}, we have
$$\sum_{i=3}^k f_n(p_i) > f_n(r) - f_n(r - \sum_{i=3}^k p_i) \ ;$$
therefore, in proving \eqref{eq:f-convexity} we can assume that $k = 2$.
Again using \eqref{eq:f-convexity-2}, we can also assume that each
nonzero~$q_j$ is equal to~$1$.
Thus, \eqref{eq:f-convexity} simplifies as follows:
\begin{equation}
\label{eq:f-convexity-simplified}
f_n(p_1) + f_n(p_2) - \ell f_n(1) - f_n(r) > 0
\,\, \text{for} \,\,
p_1 + p_2  = \ell + r, \,\, \ell < \min(p_1, p_2),
\end{equation}
which, after further simplification, reduces to showing that
\begin{equation}
\label{eq:f-convexity-simplified-2}
r 2^{r-1}  - p_1 2^{p_1-1} - p_2 2^{p_2-1} > 0
\,\, \text{for} \,\, r > \max(p_1, p_2).
\end{equation}
To finish the proof, it remains to note that
$$r 2^{r-1}  - p_1 2^{p_1-1} - p_2 2^{p_2-1} >
((r-1) 2^{r-2}  - p_1 2^{p_1-1}) +
((r-1) 2^{r-2}  - p_2 2^{p_2-1}) \geq 0 \ ,$$
as required.
\endproof

In view of Proposition~\ref{pr:geodesic-loops},
Theorem~\ref{th:normal-fan} has the following corollary.

\begin{proposition}
\label{pr:2-faces}
Every two-dimensional face of the nested polytope
is a $d$-gon for $d \in \{3, 4, 5, 6\}$.
\end{proposition}

\section{Graphical buildings}
\label{sec:graphical}

\begin{definition}[\cite{post}, Example~7.2]
\label{def:graphical-building}
Let~$\Gamma$ be a graph on the set of vertices~$S$.
Define the \emph{graphical building} $\B(\Gamma)$ as the set of all
nonempty subsets $C \subset S$ of vertices such that the induced
subgraph $\Gamma|_C$ is connected.
\end{definition}

A graphical building is indeed a building since it clearly satisfies conditions (B1) and (B2)
in Definition~\ref{def:building}.
Remembering Definition~\ref{def:restriction-contraction}, we
obtain the following proposition.

\begin{proposition}
\label{pr:restriction-contraction-graphical}
Suppose~$\B$ is a graphical building on~$S$.
\begin{enumerate}
\item
For any nonempty $C \subset S$, the restriction
$\B|_C$ is graphical.

\item
For any $C \in \B$, the contraction $C \backslash \B$
is graphical.

\end{enumerate}
\end{proposition}

\proof
Part (1) is clear since $\B|_C = \B(\Gamma|_C)$.
To show (2), it is enough to notice that
$C \backslash \B = \B(\Gamma')$, where~$\Gamma'$ is the graph
obtained from the induced subgraph $\Gamma|_{S - C}$ by adjoining
all edges $\{s, t\}$ such that each of~$s$ and~$t$ is
connected by an edge in $\Gamma$ with some vertex from~$C$.
\endproof

We now provide a characterization of
graphical buildings.

\begin{proposition}
\label{pr:graphical-characterization}
For a building~$\B$, the following conditions are equivalent:
\begin{enumerate}
\item
$\B$ is graphical.

\item
If 
$J, I_1, \dots, I_k \in \B$ are such that
$J \cup I_1 \cup \cdots \cup I_k \in \B$, then
$J \cup I_i \in \B$ for some~$i$.
\end{enumerate}
\end{proposition}

\proof
The implication $(1) \Longrightarrow (2)$ is obvious.
To prove $(2) \Longrightarrow (1)$, note that an arbitrary building~$\B$
gives rise to a graph~$\Gamma$
on the set of vertices~$S$, with $s, t \in S$ joined by an edge
whenever $\{s,t\} \in \B$.
Repeatedly using conditions (B1) and (B2) from
Definition~\ref{def:building}, we obtain the inclusion
$\B(\Gamma) \subset \B$.
It remains to show that the reverse inclusion $\B \subset \B(\Gamma)$
holds if~$\B$ satisfies~(2).

Let $C  \in \B$.
If $p = |C| \leq 2$ then $C \in \B(\Gamma)$ by the definition of~$\Gamma$;
so suppose that $p > 2$.
Repeatedly using~(2), we conclude that the elements
$s_1, \dots, s_p$ of~$C$ can be ordered in such a way that
$\{s_1, \dots, s_k\} \in \B$ for $k = 1, \dots, p$.
In particular, using induction on~$p$, we can assume that
$\{s_1, \dots, s_{p-1}\} \in \B(\Gamma)$.
Again using~(2) with $J = \{s_p\}$, and $I_1, \dots, I_k$
being the singletons $\{s_1\}, \dots, \{s_{p-1}\}$, we see that
$s_p$ is joined by an edge in~$\Gamma$ with some $s_i$ for
$1 \leq i < p$.
It follows that $C \in \B(\Gamma)$, as desired.
\endproof

We conclude with several results showing that an analogy between
the nested complexes and the cluster complexes in
\cite{yga,cfz} becomes sharper for graphical buildings.

First, as an immediate consequence of
Proposition~\ref{pr:graphical-characterization},
the nested complex of a graphical building
has the following ``clique" property.

\begin{corollary}
\label{cor:clique}
For a  graphical building~$\B$, a subset $N \subset \B - \B_{\max}$ is nested if and only if
any $I, J \in N$ satisfy \eqref{eq:B-compatibility}.
\end{corollary}

Second, for graphical buildings the property (3) of
Proposition~\ref{pr:dual-graph-exchanges} can be strengthened as
follows.

\begin{corollary}
\label{cor:exchanges-graphical}
Suppose~$\B$ is graphical, and
let $N_1$ and $N_2$ be two
maximal nested sets in $\N(\B)$ joined by an edge in the dual
graph, so that $N_1 \cap N_2 = N_1 - \{I_1\} = N_2 - \{I_2\}$.
Then, in addition to properties~(1) and~(2) of
Proposition~\ref{pr:dual-graph-exchanges}, we have:
$I_1 \cup I_2 \in (N_1 \cap N_2) \cup \B_{\max}$.
\end{corollary}

\proof
It suffices to show that, for~$\B$ graphical,
the family $\{I_3, \dots, I_k\}$ in
Proposition~\ref{pr:dual-graph-exchanges} (3) must be empty.
Suppose this is not true, and choose $k \geq 3$ to be the smallest
possible.
Then $I_i \cup I_j \notin \B$ for $3 \leq i < j \leq k$
(otherwise, replacing the pair $\{I_i, I_j\}$ with one set
$I_i \cup I_j$ would produce a smaller value of~$k$).
Applying the characterization in
Proposition~\ref{pr:graphical-characterization},
we obtain that at least one of the sets $I_1 \cup I_k$ and
$I_2 \cup I_k$ is in~$\B$.
But then if, say $I_1 \cup I_k \in \B$ then~$I_1$ and~$I_k$ cannot
belong to a nested set~$N_1$, providing the desired contradiction.
\endproof

Recall that the condition that $N_1 \cap N_2 = N_1 - \{I_1\} = N_2 - \{I_2\}$
in Corollary~\ref{cor:exchanges-graphical} means that
the corresponding cones $\RR_{\geq 0} N_1$ and $\RR_{\geq 0} N_2$ are
adjacent in the nested fan (see Corollary~\ref{cor:fan}).

\begin{corollary}
\label{cor:dependence-graphical}
In the situation of Corollary~\ref{cor:exchanges-graphical},
the linear relation \eqref{eq:dependence-concrete} takes the form
\begin{equation}
\label{eq:dependence-concrete-graphical}
\overline e_{I_1} + \overline e_{I_2} -
\sum_{J \in (\B|_{I_1 \cap I_2})_{\max}}
\overline e_{J} - \overline e_{I_1 \cup I_2} = 0 \ ;
\end{equation}
in particular, the vector $\overline e_{I_1} + \overline e_{I_2}$ belongs to
$\RR_{\geq 0} (N_1 \cap N_2)$.
\end{corollary}

Finally, for graphical buildings,
Proposition~\ref{pr:2-faces} takes the following stronger form.

\begin{proposition}[\cite{tol}, Section~2.8]
\label{pr:2-faces-graphical}
If~$\B$ is graphical then every two-dimensional face of the nested polytope
$\Pi(\B)$ is a $d$-gon for $d \in \{4, 5, 6\}$.
\end{proposition}

\proof
In view of Propositions~\ref{pr:link-decomposition}
and \ref{pr:restriction-contraction-graphical}, every
two-dimensional face of $\Pi(\B)$ corresponds to the nested
complex of a graphical building of rank~$2$.
Such complexes were listed in the proof of
Proposition~\ref{pr:geodesic-loops}.
It remains to observe that the only one among them for which the
corresponding nested polytope is a triangle (case (D1)) is
non-graphical.
\endproof

\begin{figure}[ht]
\begin{center}
\setlength{\unitlength}{1.5pt}
\begin{picture}
(120,85)(-60,-45)
\thicklines


\thinlines

\put(0,0){\circle*{1}}

\put(0,0){\vector(1,0){40}}
\put(0,0){\vector(-1,0){40}}
\put(0,0){\vector(0,1){40}}
\put(0,0){\vector(0,-1){40}}

\put(46,0){\makebox(0,0){$\overline e_1$}}
\put(-46,0){\makebox(0,0){$\overline e_3$}}
\put(0,46){\makebox(0,0){$\overline e_2$}}
\put(0,-46){\makebox(0,0){$\overline e_4$}}

\thicklines

\put(20,20){\line(0,-1){40}}
\put(20,20){\line(-1,0){40}}
\put(-20,-20){\line(1,0){40}}
\put(-20,-20){\line(0,1){40}}


\put(20,20){\circle*{2}}
\put(20,-20){\circle*{2}}
\put(-20,20){\circle*{2}}
\put(-20,-20){\circle*{2}}

\end{picture}
\qquad\qquad
\begin{picture}
(120,85)(-60,-45)
\thicklines


\thinlines

\put(0,0){\circle*{1}}

\put(0,0){\vector(1,0){40}}
\put(0,0){\vector(-1,-1){40}}
\put(0,0){\vector(0,1){40}}
\put(0,0){\vector(1,1){40}}

\put(46,0){\makebox(0,0){$\overline e_1$}}
\put(46,46){\makebox(0,0){$\overline e_{12}$}}
\put(0,46){\makebox(0,0){$\overline e_2$}}
\put(-46,-46){\makebox(0,0){$\overline e_3$}}

\thicklines

\put(30,10){\line(-1,1){20}}
\put(30,10){\line(0,-1){50}}
\put(30,-40){\line(-1,1){70}}
\put(10,30){\line(-1,0){50}}


\put(10,30){\circle*{2}}
\put(30,10){\circle*{2}}
\put(30,-40){\circle*{2}}
\put(-40,30){\circle*{2}}

\end{picture}

\end{center}
\caption{\hbox{The nested fan and polytope
in cases (D2) and (D3).}}
\label{fig:D2-D3}
\end{figure}
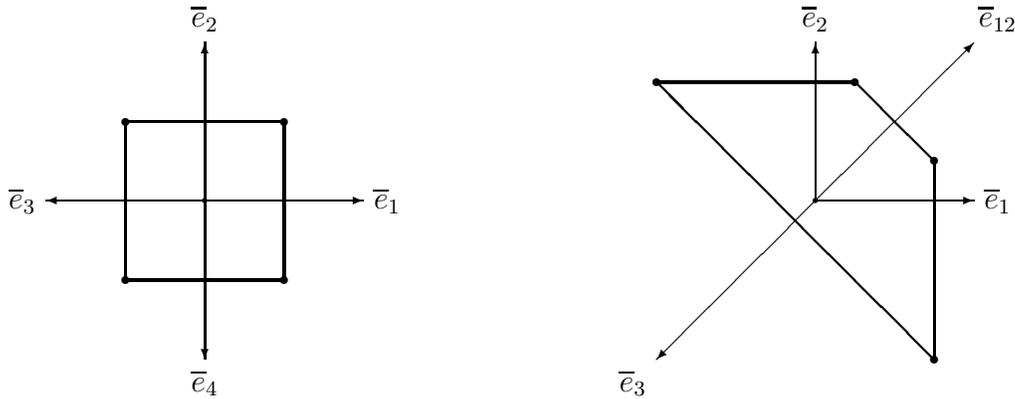

To illustrate the above results, let us compare the nested fans
and polytopes associated with rank~$2$ building sets in cases
(D2) and (D3) in the proof of Proposition~\ref{pr:geodesic-loops}.
Both cases are shown in Figure~\ref{fig:D2-D3}.
In case (D2) the building is graphical, while in (D3) it is not.
In the former case, the relations \eqref{eq:dependence-concrete-graphical}
take the form
$$\overline e_1 + \overline e_3 = 0 = \overline e_2 + \overline e_4  \ ,$$
while in the latter, we have (see \eqref{eq:dependence-concrete})
$$\overline e_1 + \overline e_2 =  \overline e_{12} = - \overline e_3 \ .$$
As a geometric consequence, in~(D2) the nested polytope is a
square, while in~(D3) it is a trapezoid.

Much more examples and pictures can be found in
\cite[Section~8]{post}.

\pagebreak[2]

\section*{Acknowledgments}
I thank Alex Postnikov and Valerio Toledano Laredo for helpful
discussions and for sharing with me the ideas and results in
their unpublished works.
Part of this work was done during my stay at the Mittag-Leffler
Institute in May - June 2005; the Institute's hospitality,
financial support ant excellent working conditions are gratefully
acknowledged.

\end{document}